\documentclass[11pt,twoside, notitlepage]{article}

\usepackage[english]{babel}
\usepackage{amsmath}
\usepackage{amsfonts}
\usepackage{amssymb}
\usepackage[all]{xy}
\usepackage{a4}
\usepackage{graphicx}
\usepackage{theorem}

\newcommand{\CC}{\mathbb{C}}

\makeatletter
\renewcommand{\subsection}{\@startsection
{subsection} {1}
{.5cm}                   
{\baselineskip} {-\fontdimen2\font plus -\fontdimen3\font minus
-\fontdimen4\font} {\normalfont\normalsize\bfseries}} \makeatother

\begin{document}
\title{Orbifolds as Groupoids: an Introduction}
\author{Ieke Moerdijk}
\date{}

\maketitle
\tableofcontents

\section*{Introduction}
The purpose of this paper is to describe orbifolds in terms of (a certain kind of) groupoids. In doing so, I hope
to convince you that the theory of (Lie) groupoids provides a most convenient language for developing the
foundations of the theory of orbifolds. Indeed, rather than defining all kinds of structures and invariants in a
seemingly ad hoc way, often in terms of local charts, one can use groupoids (which are \textit{global} objects) as
a bridge from orbifolds to \textit{standard} structures and invariants. This applies e.g. to the homotopy type of
an orbifold (via the classifying space of the groupoid), the K-theory of an orbifold (via the equivariant vector
bundles over the groupoid), the sheaf cohomology of an orbifold (via the derived category of equivariant sheaves
over the groupoid), and many other such notions. Groupoids also help to clarify the relation between orbifolds and
complexes of groups. Furthermore, the relation between orbifolds and non-commutative geometry is most naturally
explained in terms of the convolution algebra of the groupoid. In this context, and in several others, the inertia
groupoid of a given groupoid  makes its natural appearance, and helps to explain how Bredon cohomology enters the
theory of orbifolds.

Groupoids play a fundamental r\^ole in the theory of foliations, and from this point of view their use in the
context of orbifolds is only natural (cf. Theorem 3.7 below), and was exploited at an early stage, e.g. by A.
Haefliger. The more precise correspondence, between orbifolds on the one hand and proper \'etale groupoids on the
other, originates I believe with \cite{MP}.

The present paper closely follows my lecture at the Workshop in Madison, Wisconsin. In particular, the paper is of
an introductory nature, and there are no new results presented here. There are many elaborate introductions to the
theory of groupoids (e.g. \cite{CW, Mk}, and for more technical results mentioned in this paper appropriate
references are given. I have chosen to work in the context of $C^{\infty}$-manifolds, as is common in the theory
of orbifolds \cite{Sa, T}. However, it will be clear that most of the definitions and results carry over to other
categories such as the topological one. Many of the results also have analogues in algebraic geometry, in the
theory of (Deligne-Mumford and Artin) stacks.

I am most grateful to the organizers of the workshop, A. Adem, J. Morava and Y. Ruan, for inviting me to speak
there; and in particular to J. Morava for encouraging me to write up my lecture.

\section{Lie groupoids}
In this first section, we recall the basic definitions concerning
Lie groupoids.
\subsection{Groupoids.}
A groupoid is a (small) category in which each arrow is an isomorphism (\cite{CWM}, p. 20)); for example, the
category of finite sets and bijections between them. Thus, a groupoid $G$ consists of a set $G_{0}$ of
\textit{objects} and a set $G_{1}$ of \textit{arrows}, together with various structure maps. To begin with, there
are maps $s$ and $t:\xymatrix{G_{1}\ar@<.5ex>[r]\ar@<-.5ex>[r]&G_{0}}$ which assign to each arrow $g\in G_{1}$ its
\textit{source} $s(g)$ and its \textit{target} $t(g)$. For two objects $x,y\in G_{0}$, one writes
\[g:x\to y\mbox{   or   }x\overset{g}{\to}y\]
to indicate that $g\in G_{1}$ is an arrow with $s(g)=x$ and $t(g)=y$.

Next, if $g$ and $h$ are two arrows with $s(h)=t(g)$, one can form their \textit{composition} (or ``product'') $h
g$, with $s(h g)=s(g)$ and $t(h g)=t(h)$. In other words, if $g:x\to y$ and $h:y\to z$ then $h g$ is defined and
$h g: x\to z$. This composition is required to be associative. It defines a map
\begin{equation}\label{1.1.1}
G_{1}\times_{G_{0}}G_{1}\overset{m}{\to}G_{1},\qquad m(h,g)=h g
\end{equation}
on the fibered product $G_{1}\times_{G_{0}}G_{1}=\{(h,g)\in G_{1}\times G_{1}\quad|\quad s(h)=t(g)\}$.

Furthermore, for each object $x\in G_{0}$ there is a \textit{unit} (or \textit{identity}) arrow $1_{x}:x\to x$ in
$G_{1}$, which is a 2-sided unit for the composition: $g 1_{x}=g$ and $1_{x} h=h$ for any two arrows $g,h$ with
$s(g)=x=t(h)$. These unit arrows together define a map
\begin{equation}\label{1.1.2}
u:G_{0}\to G_{1},\qquad u(x)=1_{x}.
\end{equation}

Finally, for each arrow $g:x\to y$ in $G_{1}$ there exists an \textit{inverse} $g^{-1}:y\to x$, which is a 2-sided
inverse for composition: $g g^{-1}=1_{y}$ and $g^{-1} g=1_{x}$. These inverses define a map
\begin{equation}\label{1.1.3}
i:G_{1}\to G_{1},\qquad i(g)=g^{-1}.
\end{equation}

Summarizing, a groupoid consists of two sets $G_{0}$ and $G_{1}$, and five structure maps $s$, $t$, $m$, $u$, and
$i$, which are required to satisfy various identities which express that composition is associative, with a
2-sided $1_{x}$ for each $x\in G_{0}$, and a 2-sided inverse $g^{-1}$ for each $g\in G_{1}$. Before we turn to
examples, we add some smooth structure.

\subsection{Lie groupoids.}
A \textit{Lie groupoid} is a groupoid $G$ for which, in addition, $G_{0}$ and $G_{1}$ are smooth manifolds, and
the structure maps $s$, $t$, $m$, $u$ and $i$ are smooth; furthermore , $s$ and $t:
\xymatrix{G_{1}\ar@<.5ex>[r]\ar@<-.5ex>[r]&G_{0}}$ are required to be submersions (so that the domain of
definition $G_{1}\times_{G_{0}}G_{1}$ of $m$ in (\ref{1.1.1}) is a manifold). In the present context, we will
always assume $G_{0}$ and $G_{1}$ to be Hausdorff.

\subsection{Examples.}
We list some elementary standard examples.
\begin{itemize}
\item[(a)] Any groupoid can be viewed as a Lie groupoid of dimension zero; we call such a Lie groupoid \textit{discrete}.
\item[(b)] Any manifold $M$ can be viewed as a Lie groupoid in two ways: one can take $G_{1}=M=G_{0}$ which gives
a groupoid all of whose arrows are units. This is called the \textit{unit groupoid} on $M$ and denoted $u(M)$. The
other possibility, less relevant here, is to take $G_{0}=M$ and $G_{1}=M\times M$, which gives a groupoid with
\textit{exactly} one arrow $(y,x):x\to y$ from any object $x$ to any other object $y$. This is called the
\textit{pair groupoid}, and often denoted $M\times M$.
\item[(c)] Any Lie group $G$ can be viewed as a Lie groupoid, with $G_{0}=pt$ the one-point space, and
$G_{1}=G$. The composition in the groupoid is the multiplication of the group.
\item[(d)] Suppose  a Lie group $K$ acts smoothly on a manifold $M$, say from the left. One can
define a Lie groupoid $K\ltimes M$ with as objects the points $x\in M$, and as arrows $k:x\to y$ those $k$ for
which $k\cdot x=y$. Thus, $(K\ltimes M)_{0}=M$ and $(K\times M)_{1}=K\times M$, with $s: K\times M\to M$ the
projection, and $t:K\times M\to M$ the action. Composition is defined from the multiplication in the group $K$, in
the obvious way. This groupoid is called the \textit{translation groupoid} or \textit{action groupoid} associated
to the action. (For right actions there is a similar groupoid, denoted $M\rtimes K$, whose arrows $k:x\to y$ are
$k\in K$ with $x=y\cdot k$.)
\item[(e)] Let $M$ be a connected manifold. The fundamental groupoid $\Pi(M)$ of $M$ is the groupoid
with $\Pi(M)_{0}=M$ as space of objects. An arrow $x\to y$ is a homotopy class of paths from $x$ to $y$. For any
(``base'') point $x_{0}\in M$  the target $t:\Pi(M)_{1}\to M$ restricts to a map $t:s^{-1}(x_{0})\to M$ on the
space of all arrows with source $x_{0}$, and this latter map is the universal cover with base point $x_{0}$.
\end{itemize}

\subsection{Isotropy and orbits.}
Let $G$ be a Lie groupoid. For a point $x\in G_{0}$, the set of all arrows from $x$ to itself is a Lie group,
denoted $G_{x}$ and called the \textit{isotropy} or \textit{stabilizer group} at $x$. This terminology fits the
one for group actions (1.3(d)). In example 1.3(e), the stabilizer $\Pi(M)_{x}$ is the fundamental group
$\pi_{1}(M, x)$. Again inspired by group actions, the set $t s^{-1}(x)$ of targets of arrows out of $x$ is
referred to as the \textit{orbit} of $x$, and denoted $G_{x}$. This set has the structure of a smooth manifold,
for which the inclusion $G_{x}\hookrightarrow G$ is an immersion. The collection of all orbits of $G$ will be
denoted $|G|$. It is a quotient space of $G_{0}$, but in general it is not a manifold.

\subsection{Some classes of groupoids.}
There are many important properties of groupoids one can study. The following ones will be particularly relevant
here:
\begin{itemize}
\item[(a)] A Lie groupoid is called \textit{proper} if the map $(s,t):G_{1}\to G_{0}\times G_{0}$ is
a proper map \cite{B}. In a proper Lie groupoid $G$, every isotropy group is a compact Lie group.
\item[(b)] A Lie groupoid $G$ is called a \textit{foliation groupoid} if each isotropy group $G_{x}$
is \textit{discrete}. (For example, the holonomy and homotopy groupoids of a foliation have this property (see
e.g. \cite{P, CM1}.)
\item[(c)] A Lie groupoid $G$ is called \textit{\'etale} if $s$ and $t$ are local diffeomorphisms
$\xymatrix{G_{1}\ar@<.5ex>[r]\ar@<-.5ex>[r]&G_{0}}$. Any \'etale groupoid is obviously a foliation groupoid. (The
converse is true ``up to Morita equivalence'', 2.4 below.)
\item[(d)] If $G$ is an \'etale groupoid, then any arrow $g:x\to y$ in $G$ induces a well-defined
germ of a diffeomorphism $\tilde{g}:(U_{x},x)\overset{\sim}{\longrightarrow}(V_{y},y)$, as $\tilde{g}=t\circ
\hat{g}$, where $\hat{g}: U_{x}\to G_{1}$ is a section of the source map $s:G_{1}\to G_{0}$, defined on a
sufficiently small nbd $U_{x}$ of $x$ and with $\hat{g}(x)=g$. We call $G$ \textit{effective} (or
\textit{reduced}) if the assignment $g\mapsto \tilde{g}$ is faithful; or equivalently, if for each point $x\in
G_{0}$ this map $g\mapsto\tilde{g}$ defines an injective group homomorphism $G_{x}\to Diff_{x}(G_{0})$.
\end{itemize}

\section{Maps and Morita equivalence}
The main purpose of this section is to introduce the notions of (Morita) equivalence and generalized map between
Lie groupoids.

\subsection{Homomorphisms.}
Let $G$ and $H$ be Lie groupoids. A \textit{homomorphism} $\phi:H\to G$ is by definition a smooth functor. Thus, a
homomorphism consists of two smooth maps (both) denoted $\phi:H_{0}\to G_{0}$ and $\phi:H_{1}\to G_{1}$, which
together commute with all the structure maps of the groupoids $G$ and $H$. So if $h:x\to y$ is an arrow in $H$,
then $\phi(h):\phi(x)\to \phi(y)$ in $G$, etc. In the examples 1.3(b),(c) this gives the familiar notion of a map.
In example (e), any smooth map $N\to M$ is part of a unique (!) homomorphism $\Pi(N)\to \Pi(M)$. In example (d),
if $K$ acts on $M$ and $L$ on $N$, a pair consisting of a homomorphism $\alpha:K\to L$ of Lie groups and a smooth
map $f: M\to N$ which is equivariant (i.e. $f(k\cdot m)=\alpha(k)\cdot f(m)$), induces a homomorphism $K\ltimes
M\to L\ltimes N$ between translation groupoids. However, not every homomorphism is of this form.

\subsection{Natural transformations.}
Let $\phi,\psi:H\to G$ be two homomorphisms. A \textit{natural transformation} $\alpha$ \textit{from} $\phi$
\textit{to} $\psi$ (notation: $\alpha:\phi\Rightarrow \psi)$ is a smooth map $\alpha:H_{0}\to G_{1}$ giving for
each $x\in H_{0}$ an arrow $\alpha(x):\phi(x)\to \psi(x)$ in $G$, ``natural'' in $x$ in the sense that for any
$g:x\to x'$ in $G$ the identity $\psi(g)\alpha(x)=\alpha(x')\phi(g)$ should hold. We write $\phi\simeq \psi$ if
such an $\alpha$ exists. For example, if $H$ and $G$ are Lie groups (1.3(c)) $\alpha$ is just an element $g\in G$
by which the homomorphisms $\phi$ and $\psi$ are conjugate. And for example 1.3(e), if $f,g:M\to N$ are two smooth
maps and $f_{*}, g_{*}:\Pi(M)\to \Pi(N)$ are the corresponding groupoid homorphisms, then any homotopy between $f$
and $g$ defines a transformation $f_{*}\Rightarrow g_{*}$.

\subsection{Fibered products.}
Let $\phi:H\to G$ and $\psi:K\to G$ be homomorphisms of Lie groupoids. The \textit{fibered product} $H\times_{G}K$
is the Lie groupoid whose objects are triples $(y,g,z)$ where $y\in H_{0}$, $z\in K_{0}$ and $g:\phi(y)\to
\psi(z)$ in $G$. Arrows $(y,g,z)\to (y',g',z')$ in $H\times_{G}K$ are pairs $(h,k)$ of arrows, $h:y\to y'$ in $H$
and $k:z\to z'$ in $K$ with the property that $g'\phi(h)=\psi(h)g$:
\[\xymatrix{y\ar[d]_{h}&\phi(y)\ar[r]^{g}\ar[d]_{\phi(h)}&\phi(z)\ar[d]^{\phi(k)}&z\ar[d]^{k}\\
y'&\phi(y')\ar[r]_{g'}&\phi(z')&z'}\] Composition in $H\times_{G}K$ is defined in the obvious way. This fibered
product groupoid is a Lie groupoid as soon as the space
$(H\times_{G}K)_{0}=H_{0}\times_{G_{0}}G_{1}\times_{G_{0}}K_{0}$ is a manifold. (For example, this is the case
whenever the map $t\pi_{2}:H_{0}\times_{G_{0}}G_{1}\to G_{0}$ is a submersion.) There is a square of homomorphisms
\begin{equation}\label{2.1}
\xymatrix{H\times_{G} K\ar[r]^-{\pi_{2}}\ar[d]_{\pi_{1}}&K\ar[d]^{\psi}\\ H\ar[r]^{\phi}&G}
\end{equation}
which commutes up to a natural transformation, and is universal with this property.

\subsection{Equivalence.} A homomorphism $\phi:H\to G$ between Lie groupoids is called an \textit{equivalence} if
\begin{itemize}
\item[(i)] The map
\[t\pi_{1}:G_{1}\times_{G_{0}}H_{0}\to G_{0},\]
defined on the fibered product $\{(g,y)\,|\,g\in G_{1}, y\in H_{0}, s(g)=\phi(y)\}$, is a surjective submersion.
\item[(ii)] The square
\[\xymatrix{H_{1}\ar[r]^{\phi}\ar[d]_{(s,t)}&G_{1}\ar[d]^{(s,t)}\\H_{0}\times H_{0}\ar[r]^-{\phi\times \phi}&
G_{0}\times G_{0}}\] is a fibered product.
\end{itemize}

Condition (i) says in particular that every object $x$ in $G$ is connected by an arrow $g:\phi(y)\to x$ to an
object in the image of $\phi$. Condition (ii) says in particular that $\phi$ induces a diffeomorphism $H(y,z)\to
G(\phi(y),\phi(z))$, from the space of all arrows $y\to z$ in $H$ to the space of all arrows $\phi(y)\to\phi(z)$
in $G$.

Two Lie groupoids $G$ and $G'$ are said to be \textit{Morita equivalent} if there exists a third groupoid $H$ and
equivalences
\[G\overset{\phi}{\gets}H\overset{\phi'}{\to}G'.\]
This defines an equivalence relation. (To prove that the relation is transitive, one uses the fact that in a
fibered product (\ref{2.1}) above, $\pi_{1}$ is an equivalence whenever $\psi$ is.)

\subsection{Generalized maps.}
Roughly speaking, a ``generalized map'' from a Lie groupoid $H$ to a Lie groupoid $G$ is given by first replacing
$H$ by a Morita equivalent groupoid $H'$ and then mapping $H'$ into $G$ by a homomorphism of Lie groupoids. One
possible way to formalize this is via the calculus of fractions of Gabriel and Zisman \cite{GZ}. Let $\mathcal{G}$
be the category of Lie groupoids and homomorphisms. Let $\mathcal{G}_{0}$ be the quotient category, obtained by
identifying two homomorphisms $\phi$ and $\psi: \xymatrix{H\ar@<.5ex>[r]\ar@<-.5ex>[r]&G}$ iff there exists a
transformation between them (cf. 2.2). Let $W$ be the class of arrows in $\mathcal{G}_{0}$ which are represented
by equivalences (2.4) in $\mathcal{G}$. Then $W$ admits a right calculus of fractions. Let
$\mathcal{G}_{0}[W^{-1}]$ be the category of fractions, obtained from $\mathcal{G}_{0}$ by inverting all
equivalences. An arrow $H\to G$ in this category is an equivalence class of pairs of homomorphisms
\[H\overset{\varepsilon}{\gets}H'\overset{\phi'}{\to}G\]
where $\varepsilon$ is an equivalence. One may think of $H'$ as a ``cover'' of $H'$. If $\delta:H''\to H$ is a
``finer'' equivalence in the sense that there exists a homomorphism $\gamma:H''\to H'$ for which
$\varepsilon\gamma\simeq\delta$ (cf. 2.2), and if $\phi':H''\to G$ is a homomorphism or which
$\phi\gamma\simeq\phi'$, then
\[H\overset{\delta}{\gets}H''\overset{\phi'}{\to}G\]
represents the \textit{same} arrow in the category $\mathcal{G}_{0}[W^{-1}]$ of fractions. Now a generalized map
$H\to G$ is by definition an arrow in $\mathcal{G}_{0}[W^{-1}]$.

(This formulation in terms of categories of fractions  occurs in \cite{M1}. There is a more subtle ``bicategory of
fractions'' for $\mathcal{G}$ rather than $\mathcal{G}_{0}$, discussed in \cite{Pr}. There is also a way to define
these generalized morphisms in terms of principal bundles, see e.g. \cite{HS, Mr2}).

\subsection{Remark.}
Let $H$ be a Lie groupoid, and let $\mathcal{U}=\{U_{i}\}_{i\in I}$ be an open cover of $H_{0}$. Let
$H_{\mathcal{U}}$ be the Lie groupoid whose space of objects is the disjoint sum $U=\coprod U_{i}$. Write a point
in $U$ as a pair $(x,i)$ with $i\in I$ and $x\in U_{i}$. Arrows $(x,i)\to (y,j)$ in $H_{\mathcal{U}}$ are arrows
$x\to y$ in  $H$. Then the evident map $\varepsilon: H_{\mathcal{U}}\to H$ is an equivalence. Any generalized map
$H\to G$ can be represented by an open cover $\mathcal{U}$ of $H_{0}$ and the diagram
\[H\overset{\varepsilon}{\gets}H_{\mathcal{U}}\overset{\phi}{\to}G\]
for some homomorphism $\phi$. Another such representation, by $\varepsilon': H_{\mathcal{U}'}\to H$ and
$\phi':H_{\mathcal{U}'} \to G$, represents the same generalized map iff, on a common refinement of $\mathcal{U}$
and $\mathcal{U}'$, the restrictions of $\phi$ and $\phi'$ are related by a natural transformation.

\subsection{Invariance under Morita equivalence.}
For many constructions on Lie groupoids, it is important to know whether they are functorial on generalized maps.
By the universal property of the category of fractions (see \cite{GZ}), a construction is functorial on
generalized maps iff it is functorial on homomorphisms and invariant under Morita equivalence. Explicitly, let
$F:\mathcal{G}_{0}\to \mathcal{C}$ be a functor into some category $\mathcal{C}$. Then $F$ induces a functor on
the category $\mathcal{G}_{0}[W^{-1}]$ iff it sends equivalences to isomorphisms. For example, the functor mapping
a Lie groupoid $G$ to its orbit space $|G|$ has this property, since any equivalence $\phi:H\to G$ induces a
homeomorphism $|\phi|:|H|\to|G|$.

Many \textit{properties} (cf. 1.5) of Lie groupoids are invariant under Morita equivalence. For example, if $\phi:
H\to G$ is an equivalence then $H$ is proper (respectively, a foliation groupoid) iff $G$ is. Being \'etale is
\textit{not} invariant under Morita equivalence. In fact, a Lie groupoid is a foliation groupoid iff it is Morita
equivalent to an \'etale groupoid \cite{CM2}. For two Morita equivalent \'etale groupoids $G$ and $H$, one is
effective iff the other is. Thus it makes sense to call a foliation groupoid \textit{effective} (or
\textit{reduced}) iff it is Morita equivalent to a reduced \'etale groupoid.

\section{Orbifold groupoids}
In this section we explain how to view orbifolds as groupoids, and state several equivalent characterizations of
``reduced'' orbifolds.

\subsection{Definition.}
An \textit{orbifold groupoid} is a proper foliation groupoid.

For example, if $\Gamma$ is a discrete group acting properly on a manifold $M$, then $\Gamma\ltimes M$ is a proper
\'etale groupoid, hence an orbifold groupoid. Similarly, if $K$ is a compact Lie group acting on $M$, and each
stabilizer $K_{x}$ is finite, then $K\ltimes M$ is an orbifold groupoid. Observe that the slice theorem for
compact group actions gives for each point $x\in M$ a ``slice'' $V_{x}\subseteq M$ for which the action defines a
diffeomorphism $K\times_{K_{x}}V_{x} \hookrightarrow M$ onto a saturated open nbd $U_{x}$ of $x$. Then
$K_{x}\ltimes V_{x}$ is an \'etale groupoid which is Morita equivalent to $K\ltimes U_{x}$. Patching these \'etale
groupoids together for sufficiently many slices $V_{x}$ yields an \'etale groupoid Morita equivalent to $K\ltimes
M$ (cf. 2.7).

\subsection{Orbifold structures.}
If $G$ is an orbifold groupoid, its orbit space $|G|$ is a locally compact Hausdorff space. If $X$ is an arbitrary
such space, an \textit{orbifold structure} on $X$ is represented by an orbifold groupoid $G$ and a homeomorphism
$f:|G|\to X$. If $\phi: H\to G$ is an equivalence, then $|\phi|:|H|\to |G|$ is a homeomorphism and the composition
$f\circ|\phi|:|H|\to |G|\to X$ is viewed as defining an \textit{equivalent} orbifold structure on $X$. (One should
think of $G$ and $f$ as an orbifold atlas for $X$, and of $H$ and $f\circ|\phi|$ as a finer atlas.) Recall that
for any orbifold groupoid $G$, there exist equivalences $E\gets H\to G$ for which $E$ is a proper \'etale groupoid
(cf. 2.7). Thus, an orbifold structure on $X$ can always be represented by a proper \'etale groupoid $E$ and a
homeomorphism $|E|\overset{\sim}{\to}X$.

\subsection{The category of orbifolds.} An orbifold $\underline{X}$ is a space $X$ equipped with an equivalence
class of orbifold structures. A specific such structure, given by $G$ and $f:|G|\overset{\sim}{\to}X$ as in 3.2,
is then said to \textit{represent} the orbifold $\underline{X}$. For two orbifolds $\underline{X}$ and
$\underline{Y}$, represented by $(X,G,f)$ and $(Y,H,g)$ say, a \textit{map} $\underline{Y}\to \underline{X}$ is a
pair consisting of a continuous map $Y\to X$ of spaces and a generalized map $H\to G$ of orbifold groupoids, for
which the square
\[\xymatrix{|H|\ar[r]\ar[d]^{\sim}&|G|\ar[d]^{\sim}\\Y\ar[r]&X}\]
commutes. Because of the definition of generalized map, this notion of map between orbifolds is independent of the
specific representation. Notice that, by choosing a finer representation $(Y,H',g')$ for $\underline{Y}$, one can
always represent a given map $\underline{Y}\to\underline{X}$ by a continuous map $Y\to X$ together with an actual
groupoid  homomorphism $H'\to G$.

Of course, for a representation $(X,G,f)$ as above, the space $X$ is determined up to homeomorphism by the
groupoid $G$, and a similar remark applies to the maps. This means that the category of orbifolds is
\textit{equivalent} (\cite{CWM}, page 91) to a full subcategory of the category $\mathcal{G}_{0}[W^{-1}]$ of Lie
groupoids and generalized maps, viz. the category determined by the orbifold groupoids. It is often easier to work
explicitly with this category. Notice also that a fibered product $H\times_{G}K$ in (4) is an orbifold groupoid
whenever $H$, $G$ and $K$ are (and the transversality conditions are met for $H\times_{G}K$ to be smooth.) This
defines a notion of fibered product of orbifolds.

\subsection{Local charts.}
Let $G$ be a Lie groupoid. For an open set $U\subseteq G_{0}$, we write $G|_{U}$ for the full subgroupoid of $G$
with $U$ as a space of objects. In other words, $(G|_{U})_{0}=U$ and $(G|_{U})_{1}=\{g:x\to y\,|\,g\in G_{1}\mbox{
and }x,y\in U\}$. If $G$ is proper and \'etale, then for each $x\in G_{0}$ there exist arbitrary small
neighborhoods $U$ of $x$ for which $G|_{U}$ is isomorphic to $G_{x}\ltimes U$ for an action of the isotropy group
$G_{x}$ on the neighborhood $U$ (see e.g. \cite{MP}). In particular, such a $U$ determines an open set
$|U|\subseteq |G|$, for which $|U|$ is just the quotient of $U$ by the action of the finite group $G_{x}$. Note
also that in the present smooth context, one can choose coordinates so that $U_{x}$ is a Euclidean ball on which
$G_{x}$ acts linearly. In the literature, one often \textit{defines} an orbifold in term of local quotients. It is
also possible to describe maps between orbifolds or orbifold groupoids (3.3) in terms of these local charts. This
leads to the notion of \textit{strict map} of \cite{MP}, or the equivalent notion of \textit{good map} of Chen and
Ruan.

\subsection{Embedding categories.}
Let $G$ be an \'etale groupoid, and let $\mathcal{B}$ be a basis for $G_{0}$. The embedding category
$\underline{Emb}_{\mathcal{B}}(G)$ of $G$ determined by $\mathcal{B}$ is a discrete (small) category with elements
of $\mathcal{B}$ as objects. For $U,V\in\mathcal{B}$, an arrow $U\to V$ in the embedding category is a smooth map
$\sigma: U\to G_{1}$ with the property that $\sigma(x)$ is an arrow from $x$ to some point $y=t \sigma(x)$ in $V$,
and such that the corresponding map $t\circ\sigma: U\to V$, sending $x$ to $y$, is an embedding of $U$ into $V$.
For two such arrows $\sigma:U\to V$ and $\tau:V\to W$, their composition $\tau\circ\sigma:U\to W$ is defined in
the obvious way from the composition in $G$, by
\[(\tau\circ\sigma)(x)=\tau(t\sigma(x))\sigma(x).\]
It is possible to reconstruct the groupoid $G$ from the manifold $G_{0}$, its basis $\mathcal{B}$, and this
embedding category $\underline{Emb}_{\mathcal{B}}(G)$. Such embedding categories play an important r\^{o}le, e.g.
in describing the homotopy type of $G$ and its characteristic classes; see also \S4. It should also be observed
that it is often easier to construct the embedding category $\underline{Emb}_{\mathcal{B}}(G)$ directly, and not
in terms of the groupoid $G$. This applies in particular to situations where an orbifold is given by local charts
and ``embeddings'' between them, as in \cite{Sa}.

\subsection{Reduced orbifolds.}
An orbifold structure is said to be \textit{reduced} if it is given by an \textit{effective} orbifold groupoid $G$
(cf. 1.3(d)). An orbifold is called \textit{reduced} if it is represented by a reduced orbifold structure. Any
reduced orbifold can be represented by an effective proper \'etale groupoid $G$. If $G$ is such a groupoid, then
the local charts $G_{x}\ltimes U$ of 3.4 have the property that $G_{x}$ acts effectively on $U$. Notice that if
$G$ is an arbitrary \'etale groupoid, then there is an evident quotient $G\twoheadrightarrow \tilde{G}$, which is
an effective \'etale groupoid. The space $\tilde{G}_{0}$ of objects is the same as that of $G$, i.e.
$\tilde{G}_{0}=G_{0}$, and $|G|=|\tilde{G}|$. Thus, an arbitrary orbifold structure on a space $X$ always has an
associated reduced structure on the same space $X$.

The following theorem lists some of the ways in which reduced orbifolds occur in nature.

\subsection{Theorem.}
\textit{For a Lie groupoid $G$, the following conditions are equivalent.
\begin{itemize}
\item[(i)] $G$ is Morita equivalent to an effective proper \'etale groupoid.
\item[(ii)] $G$ is Morita equivalent to the holonomy groupoid of a foliation with compact leaves
and finite holonomy.
\item[(iii)] $G$ is Morita equivalent to a translation groupoid $K\ltimes M$, where $K$ is a compact
Lie group acting on $M$, with the property that each stabilizer $K_{x}$ is finite and acts effectively on the
normal bundle.
\end{itemize}}

More explicitly, in (iii) the stabilizer $K_{x}$ acts on the tangent space $T_{x}(M)$, and this action induces an
action on the quotient bundle $N_{x}$ of tangent vectors normal to the orbits of the action. This latter action
must be effective in (iii).

This theorem is stated and proved  in some detail, e.g. in \cite{MP}. However, each of the equivalences seems to
be a folklore result. The equivalence between (i) and (iii) is proved by the ``frame bundle trick'', while the
equivalence between (i) and (ii) is closely related to the Reeb stability theorem for foliations.

\subsection{Complexes of groups.}
There is a variation on the construction of embedding categories $\underline{Emb}_{\mathcal{B}}(G)$ which works
well for reduced orbifolds, and is related to the theory of complexes of groups \cite{Hf, Bs, Sr}. In general, a
complex of (finite) groups over a small category $I$ is a pseudofunctor $F$ from $I$ into the category of (finite)
groups and injective homomorphisms. This means that $F$ assigns to each object $i$ of $I$ a group $F(i)$, and to
each arrow $\sigma:i\to j$ in $I$ an injective homomorphism $F(\sigma):F(i)\to F(j)$. Furthermore, this assignment
is ``pseudo-functorial'' in the sense that for any composable pair $\sigma:i\to j$ and $\tau:j\to k$ of arrows in
$I$, there is an element $g_{\tau,\sigma}\in F(k)$ for which
\[\hat{g}_{\tau,\sigma} F(\tau\sigma)=F(\tau)F(\sigma):F(i)\to F(k),\]
where $\hat{g}_{\tau,\sigma}:F(k)\to F(k)$ is the inner automorphism $\hat{g}_{\tau,\sigma}(x)=g_{\tau,\sigma}x
g_{\tau,\sigma}^{-1}$. Moreover, these $\hat{g}_{\tau,\sigma}$ should satisfy a coherence condition: for
$i\overset{\sigma}{\to}j\overset{\tau}{\to}k\overset{\rho}{\to}l$ in $I$,
\begin{equation}\label{3.8.1}
F(\rho)(g_{\tau,\sigma})g_{\rho,\tau\sigma}=g_{\rho,\tau}g_{\rho\tau,\sigma}
\end{equation}
(an identity in the group $F(l)$). We will assume that the pseudo functor $F$ is ``normal'' in the sense that for
any identity arrow $1:i\to i$ in $I$, it takes the value $id:F(i)\to F(i)$. To such a pseudofunctor $F$ on $I$ one
can associate a total category $\int_{I}F$ by the well known ``Grothendieck construction''. In this particular
case, this is a category with the same objects as $I$. Arrows $i\to j$ in $\int_{I}F$ are pairs $(\sigma,g)$ where
$\sigma:i\to j$ and $g\in F(j)$. The composition of $(\sigma,g):i\to j$ and $(\tau,h):j\to k$ in $\int_{I}F$ is
defined as $(\tau\sigma, h F(\tau)(g) g_{\tau,\sigma})$; this is an associative operation, by (\ref{3.8.1}) above.

\subsection{The complex of groups associated to the reduced orbifold.}
Let $G$ be a reduced \'etale groupoid, with orbit space $|G|$. By 3.4, there is a basis $\mathcal{B}$ for $G_{0}$,
consisting of simply connected open sets $U\subseteq G_{0}$, such that for each $U$ in $\mathcal{B}$ the
restricted groupoid $G|_{U}$ is isomorphic to a translation groupoid $G_{U}\ltimes U$, where $G_{U}=G_{x}$ is the
stabilizer group at a suitable point $x\in U$. For such a basis $\mathcal{B}$, the images $q(U)\subseteq|G|$ under
the quotient map $q:G_{0}\to |G|$ form a basis $\mathcal{A}$ for $|G|$. Fix such a basis $\mathcal{A}$, and choose
for each $A\in \mathcal{A}$ a specific $U=\tilde{A}\in\mathcal{B}$ for which $A=q(\tilde{A})$. View the partially
ordered set $\mathcal{A}$ as a category, with just one arrow $A\to A'$ iff $A\subseteq A'$, as usual. If
$A\subseteq A'$, then the fact that $\tilde{A}$ is simply connected  readily implies that there exists an arrow
$\lambda_{A,A'}:\tilde{A}\to\tilde{A}'$ in the embedding category $\underline{Emb}_{\mathcal{B}}(G)$. This arrow
induces a group homomorphism
\[(\lambda_{A,A'})_{*}:G_{\tilde{A}}\to G_{\tilde{A'}},\]
of local groups belonging to the ``charts'' $\tilde{A}$ and $\tilde{A}'$ in $\mathcal{B}$. If $A\subseteq
A'\subseteq A''$, then $\lambda_{A',A''}\circ\lambda_{A,A'}$ and $\lambda_{A,A''}$ differ by conjugation by an
element $g$ of $G_{A''}$, and this element is unique if $G$ is effective. Thus, the assignment $A\mapsto
G_{\tilde{A}}$ and $(A\subseteq A')\mapsto (\lambda_{A,A'})_{*}$ is part of a uniquely determined complex of
groups, which we denote by $G_{\mathcal{A}}$. Its total category $\int_{\mathcal{A}}G_{\mathcal{A}}$ is smaller
than, but in fact categorically equivalent to, the embedding category $\underline{Emb}_{\mathcal{B}}(G)$.

\section{The classifying space}
The purpose of this section is to introduce the homotopy type, as the classifying space of a representing
groupoid.

\subsection{The nerve of a groupoid.}
Let $G$ be a Lie groupoid, given by manifolds $G_{0}$ and $G_{1}$ and the various structure maps. Let $G_{n}$ be
the iterated fibered product
\[G_{n}=\{(g_{1},\ldots,g_{n})\,|\,g_{i}\in G_{1},s(g_{i})=t(g_{i+1})\mbox{ for }i=1,\ldots,n-1\}.\]
In other words, $G_{n}$ is the manifolds of composable strings
\[x_{0}\overset{g_{1}}{\gets}x_{1}\gets\ldots\overset{g_{n}}{\gets}x_{n}\]
of arrows in $G$. These $G_{n}$ together have the structure of a simplicial manifold, called the \textit{nerve} of
$G$. The face operator $d_{i}:G_{n}\to G_{n-1}$ for $i=0,\ldots,n$ are given by ``$d_{i}=\,\mbox{delete }x_{i}$'';
so $d_{0}(g_{1},\ldots,g_{n})=(g_{2},\ldots,g_{n})$, and $d_{n}(g_{1},\ldots,g_{n})=(g_{1},\ldots,g_{n-1})$, while
$d_{i}(g_{1},\ldots,g_{n})=(g_{1},\ldots, g_{i}g_{i+1},\ldots,g_{n})$ for $0<i<n$.

The same definition of course makes sense for groupoids in categories others than that of differentiable
manifolds. In particular,  if $G$ is only a topological groupoid then the $G_{n}$ form a simplicial topological
space. The definition doesn't involve the inverse of a groupoid either, so makes sense for topological categories,
for example.

\subsection{The classifying space.} (\cite{Sg2}) For a simplicial space $X_{\bullet}$, we write
$|X_{\bullet}|$ for its geometric realization. This is a space obtained by gluing the space $X_{n}\times
\Delta^{n}$ along the simplicial operators; here $\Delta^{n}$ is the standard $n$-simplex. If $X_{\cdot}$ is a
simplicial \textit{set} (i.e. each $X_{n}$ has the discrete topology) then $|X_{\bullet}|$ is a CW-complex.
However in general this is not the case and there are various subtleties involved in this definition. In
particular, if the degeneracies $X_{n-1}\hookrightarrow X_{n}$ aren't cofibrations, one should take the ``thick''
realization, cf. \cite{Sg2}. For a Lie groupoid $G$, its classifying space $B G$ is defined as the geometric
realization of its nerve,
\[B G=|G_{\bullet}|\]
This notation should not be confused with the notation $|G|$ for the orbit space of $G$.

\subsection{Homotopy type of an orbifold.}
An important basic property of the classifying space construction is that an equivalence $\phi:H\to G$ (cf. 2.4)
induces a weak homotopy equivalence $B \phi:B H\to B G$. This means that for any point $y\in H_{0}$ this
equivalence $H\to G$ induces an isomorphism of homotopy groups $\pi_{n}(B H,y)\to \pi_{n}(B G,\phi y)$. Thus, if
$\underline{X}$ is an orbifold, one can define its \textit{homotopy type} as that of $B G$ where $G$ is any
orbifold structure representing $\underline{X}$; the definition
\[\pi_{n}(\underline{X},x)=\pi_{n}(B G,\tilde{x})\]
is then independent of the groupoid $G$ and the base point $x\in X$, and of a ``lift'' $\tilde{x}\in G_{0}$ for
which $q(\tilde{x})\in |G|$ is mapped to $x\in X$ by the given homeomorphism $|G|\to X$.

\subsection{Local form of the classifying space.}
Let $G$ be a proper \'etale groupoid representing an orbifold $\underline{X}$. The map $q:G_{0}\to |G|$ to the
orbit space induces an obvious map, still denoted $q:B G\to |G|$. By 3.4 (and as in 3.9) there is an open cover
$|G|=\bigcup V_{i}$ such that $G|_{q^{-1}(V_{i})}$ is Morita equivalent to $G_{i}\ltimes U_{i}$, where $G_{i}$ is
a finite group acting on an open set $U_{i}\subseteq G_{0}$. One can choose the $V_{i}$ and $U_{i}$ to be
contractible. Thus $q^{-1}(V_{i})\subseteq B G$ is homotopy equivalent to $B(G_{i}\ltimes U_{i})$, which is the
Borel space $E G_{i}\times_{G_{i}}U_{i}$, homotopy equivalent to $B G_{i}$ if we choose $U_{i}$ to be
contractible. Thus, over $V_{i}$ the map $q: B G\to |G|$ restricts to $E G_{i}\times_{G_{i}} U_{i}\to
U_{i}/G_{i}\cong V_{i}$. In particular, it follows, for example, that $q: BG\to |G|$ induces isomorphisms in
rational cohomology.

\subsection{Other models for the classifying space.}
Let $G$ be an \'etale groupoid. Let $\mathcal{B}$ be a basis for $G_{0}$ consisting of contractible open sets, and
let $\underline{Emb}_{\mathcal{B}}(G)$ be the associated embedding category (3.5). In \cite{M3} it is proved that
there is a weak homotopy equivalence
\[B G\simeq B \underline{Emb}_{\mathcal{B}}(G).\]
Notice that $\underline{Emb}_{\mathcal{B}}(G)$ is a discrete small category, so the right-hand side is a
CW-complex.

If $G$ is \textit{proper}, one can choose the basis $\mathcal{B}$ to consist of linear charts as in 3.4. Then both
$B$ and $q(B)$ are contractible. Let $\mathcal{A}=\{q(B)\,|\,B\in \mathcal{B}\}$ as in 3.9. Then $\mathcal{A}$ is
a basis for $|G|$, and when we view $\mathcal{A}$ as a category, $B \mathcal{A}\simeq |G|$. If $G$ is moreover
\textit{effective}, then there is an associated complex of groups $G_{\mathcal{A}}$. The classifying space of its
total category $\int_{\mathcal{A}}G_{\mathcal{A}}$ is again a model for the homotopy type of the orbifold, because
there is a weak homotopy equivalence
\[B G\simeq B(\int_{\mathcal{A}}G_{\mathcal{A}}).\]
The quotient map $q:BG\to |G|$ is modelled by the projection functor $\int_{\mathcal{A}}G_{\mathcal{A}}
\to\mathcal{A}$, in the sense that there is a diagram
\[\xymatrix{B(\int_{\mathcal{A}}G_{\mathcal{A}})\ar[r]^{\sim}\ar[d]&B G\ar[d]^{q}\\ B \mathcal{A}\ar[r]^{\sim}&|G|}\]
which commutes up to homotopy.

\section{Structures over orbifolds}
In this section we show how the language of groupoids leads to a uniform definition of structures ``over''
orbifolds, like covering spaces, vector bundles, principal bundles, sheaves, etc. We continue to work in the
smooth context, although everything in this section obviously applies equally to topological groupoids.

\subsection{$G$-spaces.}
Let $G$ be a Lie groupoid. A (right) \textit{$G$-space} is a manifold $E$ equipped with an action by $G$. Such an
action is given by two maps, $\pi:E\to G_{0}$ and $\mu:E\times_{G_{0}}G_{1}\to E$. The latter map is defined on
pairs $(e,g)$ with $\pi(e)=t(g)$, and written $\mu(e,g)=e\cdot g$. It satisfies the usual identities for an
action, viz. $\pi(e\cdot g)=s(g)$, $e\cdot 1_{x}=e$ and $(e\cdot g)\cdot h=e\cdot(g h)$ for
$z\overset{h}{\to}y\overset{g}{\to}x$ in $G$ and $e\in E$ with $\pi(e)=x$. For two such $G$-spaces $E=(E,\pi,\mu)$
and $E'=(E',\pi',\mu')$, a map of $G$-spaces $\alpha:E\to E'$ is a smooth map which commutes with the structure,
i.e. $\pi'\alpha=\pi$ and $\alpha(e\cdot g) =\alpha(e)\cdot g$. This defines a category \textit{($G$-spaces)}.

If $\phi:H\to G$ is a homomorphism of groupoids, there is an obvious functor
\begin{equation}\label{5.1.1}
\phi^{*}: \mbox{\textit{($G$-spaces)}}\to \mbox{\textit{($H$-spaces)}}
\end{equation}
mapping $E$ to the pullback $E\times_{G_{0}}H_{0}$ with the induced action. If $\phi$ is an equivalence (2.4) then
this functor $\phi^{*}$ is an equivalence of categories. Thus, up to equivalence of categories, the category
\textit{($G$-spaces)} only depends on the Morita equivalence class of $G$.

\subsection{$G$-spaces as groupoids.}
If $E$ is a $G$-space, then one can form the translation groupoid $E\rtimes G$ whose objects are point in $E$, and
whose arrows $g:e'\to e$ are arrows $g:\pi(e')\to\pi(e)$ in $G$ with $e\cdot g=e'$. In other words, $(E\rtimes
G)_{0}=E$ and $(E\rtimes G)_{1}=E\times_{G_{0}}G_{1}$, while the source and target of $E\rtimes G$ are the action
$\mu$ and the projection $E\times_{G_{0}}G_{1}\to E$. There is an obvious homomorphism of groupoids
$\pi_{E}:E\rtimes G\to G$. Observe also that for a homomorphism $\phi:H\to G$, the square
\[\xymatrix{\phi^{*}(E)\rtimes H\ar[r]\ar[d]&E\ar[d]\\H\ar[r]&G}\]
is a fibered product (2.3) up to Morita equivalence. Notice that, while at the groupoid level the fiber of
$E\rtimes G\to G$ over $x\in G_{0}$ is the fiber $\pi^{-1}(x)$, at the level of orbit spaces $|E\rtimes G|\to |G|$
the fiber is $\pi^{-1}(x)/G_{x}$. It is easy to see that
\begin{itemize}
\item[(i)]
If $G$ is \'etale then so is $E\rtimes G$.
\item[(ii)]
If $G$ is a foliation groupoid then so is $E\rtimes G$.
\item[(iii)]
If $E$ is Hausdorff and $G$ is proper then $E\rtimes G$ is proper.
\end{itemize}
In particular, if $G$ represents an orbifold $\underline{X}$ then any Hausdorff $G$-space $E$ represent an
orbifold $\underline{E}\to \underline{X}$ over $\underline{X}$.

We now mention some examples:

\subsection{Covering spaces.}
A covering space over $G$ is a $G$-space $E$ for which $\pi:E\to G_{0}$ is a covering projection. The full
subcategory of \textit{($G$-spaces)} consisting of covering spaces is denoted $Cov(G)$. For a groupoid
homomorphism $\phi:H\to G$, the functor $\phi^{*}:H\to G$ in (\ref{5.1.1}) restricts to a functor
\[\phi^{*}:Cov(G)\to Cov(H)\]
and this is an equivalence of categories whenever $\phi$ is an equivalence. So, up to equivalence of categories,
there is a well defined category $Cov(\underline{X})$ of covering spaces of an orbifold $\underline{X}$.

If $G$ is proper and \'etale while $E$ is a covering space of $G$, then $E\rtimes G$ is again proper and \'etale.
The map $|E\rtimes G|\to |G|$ of orbit spaces has local charts (cf. 3.4) $\tilde{U}/\tilde{\Gamma}\to U/\Gamma$
where $\tilde{U}\cong U$, $\Gamma=G_{x}$ is a finite group, and $\tilde{\Gamma}\subseteq \Gamma$ is a subgroup.
This explains the relation to the definition of covering spaces of orbifolds given in \cite{T}.

Suppose $|G|$ is connected. If $x\in G_{0}$ is a base point, there is a fiber functor $F_{x}:Cov(G)\to Sets$,
mapping $E$ to $E_{x}=\pi^{-1}(x)$. Grothendieck's Galois theory applies and gives a unique (up to isomorphism)
group $\pi_{1}(G, x)$ for which there is an equivalence of categories
\[Cov(G)\overset{\sim}{\to}\pi_{1}\mbox{\textit{$(G,x)$-sets}}\]
by which $F_{x}$ corresponds to the forgetful functor $\pi_{1}\mbox{\textit{$(G,x)$-sets}}\to
\mbox{\textit{sets}}$. It follows from earlier remarks that if $\phi:H\to G$ is an equivalence and $y\in H_{0}$
then $\phi$ induces an isomorphism $\pi_{1}(H,y)\to\pi_{1}(G,\phi(y))$. In particular, the fundamental group
$\pi_{1}(G,x)$ only depends on the Morita equivalence class of $G$.

It is also possible to describe the fundamental group in terms of ``paths'', i.e. generalized maps $[0,1]\to G$ as
is done in work of Haefliger and of Mrcun. Alternatively it is not difficult to see that the map $E\mapsto
B(E\rtimes G)$ induces an equivalence of categories between covering spaces of $G$ and covering spaces (in the
usual sense) of $B G$. It follows that $\pi_{1}(G)=\pi_{1}(B G)$. (This is true more generally for arbitrary
\'etale groupoids, see \cite{M4}).

\subsection{Vector bundles.}
A \textit{vector bundle} over $G$ is a $G$-space $E$ for which $\pi:E\to G_{0}$ is a vector bundle, and the action
of $G$ on $E$ is fiberwise linear. (In particular, each fiber $E_{x}$ is a linear representation of the stabilizer
$G_{x}$.) Write $Vect(G)$ for the category of vector bundles over $G$. Again, if $G$ is Morita equivalent to $H$
then $Vect(G)$ is equivalent to $Vect(H)$, so up to equivalence of categories there is a well defined category of
vector bundles over an orbifold $\underline{X}$, denoted $Vect(\underline{X})$. From $Vect(\underline{X})$ one can
build the Grothendieck group $K(\underline{X})$ in the usual way, leading to $K$-theory of orbifolds; see e.g.
\cite{AR}. If $\underline{X}$ is reduced, then by Theorem 3.7 $Vect(\underline{X})$ is equivalent to the category
of $K$-equivariant vector bundles on a manifold $M$, and we are back to equivariant $K$-theory for compact Lie
groups, \cite{Sg3}.

\subsection{Principal bundles.}
Let $L$ be a Lie group. A principal $L$-bundle over $G$ is a $G$-space $P$ with a left action $L\times P\to P$
which makes $\pi:P\to G_{0}$ into a principal $L$-bundle over the manifold $G_{0}$, and is compatible with the
$G$-action in the sense that $(l\cdot p)\cdot g=l\cdot(p\cdot g)$, for any $p\in P$, $l\in L$ and $g:x\to y$ with
$y=\pi(p)$. For such principal bundles, one can construct characteristic classes exactly as for manifolds, via a
``Chern-Weil'' map from the algebra $Inv(\mathfrak{l})$ of invariant polynomials on the Lie algebra $\mathfrak{l}$
of $L$ into $H^{*}(B G)$. For a recent geometric description of these in the general context of \'etale groupoids,
see \cite{CM2}.

\subsection{Sheaves of sets and topoi.}
A $G$-space $E$ is called a sheaf of sets if $\pi:E\to G_{0}$ is a local diffeomorphism. The category of all
sheaves of sets is denoted $Sh(G)$. This category is a \textit{topos}, from which one can recover the groupoid $G$
as a topological groupoid up to Morita equivalence. More precisely, if $G$ is \'etale, the sheaf of germs of
smooth functions $G_{0}$ has the structure of a $G$-sheaf, denoted $\mathcal{A}_{G}$. Thus
$(Sh(G),\mathcal{A}_{G})$ is a ringed topos. Up to Morita equivalence one can reconstruct $G$ with its smooth
structure from this ringed topos. More generally, morphisms of ringed topoi
$(Sh(H),\mathcal{A}_{H})\to(Sh(G),\mathcal{A}_{G})$ correspond \textit{exactly} to generalized maps of groupoids
$H\to G$. Thus, the category of orbifolds could also have been introduced as a full subcategory of the category of
ringed topoi. In the language of \cite{SGA4}, a ringed topos represents an orbifold iff it is a ``separated smooth
\'etendue''. This is the viewpoint taken in \cite{MP}.

A topos is defined more generally as the category of sheaves on a \textit{site}. It is easy to give an explicit
site for ``orbifold topoi'' of the form $Sh(G)$. Indeed, a site for $Sh(G)$ is formed by any embedding category
$Emb_{\mathcal{B}}(G)$ with its evident Grothendieck topology, where a family $\{\sigma_{i}:U_{i}\to V\}$ covers
iff $V=\bigcup \sigma_{i}(U_{i})$. These sites are useful in the description of ``higher structures'' over
orbifolds such as gerbes \cite{LU} and 2-gerbes. They are also relevant for the definition of the ``\'etale
homotopy groups'' \cite{AM} of the topos $Sh(G)$. By the general comparison theorem of \cite{M4}, these \'etale
homotopy groups coincide with the homotopy groups of the classifying space $B G$. For yet another approach to the
homotopy theory of orbifolds, see \cite{Ch}.

\section{Cohomology and inertia orbifolds}
In this section we will briefly discuss sheaf cohomology of orbifolds, and its relation to the cohomology of (the
various models for) its classifying space. By way of example, we will in particular discuss the case of the
``inertia orbifold'' of a given orbifold, and its relation to K-theory and non-commutative geometry. We begin by
mentioning, as a continuation of the previous section, yet another kind of structure over an orbifold.

\subsection{Abelian sheaves.}
An abelian sheaf over an orbifold groupoid is a $G$-sheaf $A$ (cf 5.6) for which the fibers $A_{x}=\pi^{-1}(x)$ of
$\pi:A\to G_{0}$ have the structure of an abelian group, and this group structure varies continuously in $x$ and
is preserved by the action of arrows in $G$. For example, germs of local sections of a vector bundle over $G$ form
an abelian sheaf. The category $Ab(G)$ of these abelian sheaves is a nice abelian category with enough injectives,
and up to equivalence of categories it only depends on the Morita equivalence class of $G$. In fact, $Ab(G)$ is
the category of abelian group objects in the topos $Sh(G)$, so enjoys all the general properties of the
\cite{SGA4} framework. In particular, the derived category of complexes of abelian $G$-sheaves is a convenient
model for the derived category $D(\underline{X})$ of an orbifold, in terms of which one can define all the usual
operations, prove change-of-base formulas, etc. We will not discuss these matters in detail here, but instead
refer to \cite{M2} and references cited here. Some of the following remarks 6.2-6.3 are just special cases of this
general framework.

\subsection{Cohomology.}
Let $G$ be an orbifold groupoid, and let $A$ be an abelian $G$-sheaf. The cohomology $H^{n}(G,A)$ for $n\geq 0$ is
defined as the cohomology of the complex $\Gamma_{inv}(G,I^{\bullet})$. Here $0\to A\to I^{0} \to I^{1}\to\ldots$
is a resolution by injective abelian $G$-sheaves; furthermore, for an arbitrary $G$-sheaf $B$, $\Gamma_{inv}(G,B)$
denotes the group of sections $\beta:G_{0}\to B$ of $\pi:B\to G_{0}$ which are invariant under the $G$-action,
i.e. $\beta(y)\cdot g=\beta(x)$ for any arrow $g:x\to y$ in $G$. These cohomology groups are functorial in $G$, in
the sense that any homomorphism of groupoids $\phi:G'\to G$ induces group homomorphisms $\phi^{*}:H^{n}(G,A)\to
H^{n}(G',\phi^{*}A)$. If $\phi$ is an equivalence, then this map is an isomorphism. Thus if $\underline{X}$ is an
orbifold, there are well defined cohomology groups $H^{n}(\underline{X},A)$, constructed in terms of any
representing groupoid. The cohomology groups $H^{n}(G,A)$ are also closely related to the (sheaf) cohomology
groups of the classifying space $B G$, and to the cohomology groups of the small category
$\underline{Emb}_{\mathcal{B}}(G)$ (cf. 3.5).

\subsection{Direct image functor.}
Let $\phi:H\to G$ be a homomorphism between orbifold groupoids, which we may as well assume to be \'etale. Any
$H$-sheaf $B$ defines a $G$-sheaf $\phi_{*}(B)$, whose sections over an open $U\subseteq G_{0}$ are given by the
formula
\[\Gamma(U,\phi_{*}B)=\Gamma_{inv}(U/\phi,B).\]
Here $U/\phi$ is the ``comma'' groupoid whose objects are pairs $(y,g)$ with $y\in H_{0}$ and $g:x\to\phi(y)$ an
arrow in $G$. Arrows $(y,g)\to (y',g')$, where $g':x'\to\phi(y')$, in this groupoid $U/\phi$ only exist if $x=x'$
and arrows are $h:y\to y'$ in $H$ with $\phi(h)g=g'$ in $G$. This groupoid $U/\phi$ is again a proper \'etale
groupoid, and $B$ can be viewed as a $U/\phi$-sheaf, via the evident homomorphism $U/\phi\to H$. The higher
derived functors $\phi_{*}$ can be described by the formula
\[R^{i}\phi_{*}(B)_{x}=\underset{\to}{lim}_{x\in U}H^{i}(U/\phi,B),\]
for the stalk at an arbitrary point $x\in G_{0}$.

Let us mention some examples:

(a) Let $S$ be a $G$-space and consider the projection homomorphism $\pi:S\rtimes G\to G$. Then $U/\pi$ is Morita
equivalent to the space $\pi^{-1}(U)\subseteq S$ (viewed as a unit groupoid, see 1.3(b)), and $R^{i}\pi_{*}(B)$ is
the usual right derived functor along $S\to G_{0}$; or more precisely, the diagram
\[\xymatrix{S\ar[d]_{\pi}\ar[r]^-{v}&S\rtimes G\ar[d]^{\pi}\\G_{0}\ar[r]^{u}&G}\]
with unit groupoids on the left, has the property that
\[u^{*}\circ R \pi_{*}=R \pi_{*}\circ v^{*}.\]

(b) Let $|G|$ be the orbit space of $G$, as before. In general $|G|$ is not a manifold, but we can view it as a
topological unit groupoid and consider the quotient homomorphism of groupoids $q:G\to |G|$. Then for any $G$-sheaf
$A$, and any $x\in G_{0}$, there is a natural isomorphism
\[R^{i}q_{*}(A)_{q(x)}=H^{i}(G_{x},A_{x}).\]
On the right, we find the cohomology of the finite group $G_{x}$ with coefficients in the stalk $A_{x}$. In
particular, $R^{i} q_{*}=0$ for $i>0$ if the coefficient sheaf $A$ is a sheaf of $\mathbb{Q}$-modules.

(c) Let $\phi:H\to G$ as above. Call $\phi$ \textit{proper} if $|G_{0}/\phi|\to G_{0}$ is a proper map. For such a
$\phi$ we have
\[R^{i} \phi_{*}(B)_{x}=H^{i}(x/\phi,B).\]
Here $x/\phi\subseteq U/\phi$ is the comma groupoid with objects $(x\overset{g}{\to}\phi(y))$ and arrows as in
$U/\phi$. The $H$-sheaf $B$ is viewed as an $(x/\phi)$-sheaf via the pullback along the projection $x/\phi\to H$.

\subsection{The inertia orbifold.}
Let $G$ be an orbifold groupoid, and consider the pullback of spaces
\[\xymatrix{S_{G}\ar[r]\ar[d]_{\beta}&G_{1}\ar[d]^{(s,t)}\\G_{0}\ar[r]^-{diag}&G_{0}\times G_{0}.}\]
Thus, $S_{G}=\{g\in G_{1}: s(g)=t(g)\}$ is the space of ``loops'' in $G$. The map $\beta$ sends such a loop
$g:x\to x$ to its ``base point'' $\beta(g)=x$. This map $\beta$ is proper, since $(s,t)$ is proper. Moreover, the
space $S_{G}$ is in fact a manifold. This is clear from the local charts (3.4) in the case $G$ is \'etale, and it
follows for general orbifold groupoids by considering local charts as in \cite{CM1}. The space $S_{G}$ has the
structure of a $G$-space, with the action defined by conjugation. Let
\[\Lambda(G)=S_{G}\rtimes G.\]
Then $\Lambda(G)$ is an orbifold groupoid, and $\beta$ induces a proper (cf. 6.3(c)) homomorphism
\[\beta:\Lambda(G)\to G.\]
If $\phi:H\to G$ is an equivalence then $\phi^{*}(S_{G})=S_{H}$, and $\phi$ induces an equivalence
$\phi^{*}(S_{G})\rtimes H\to S_{G}\rtimes G$, hence an equivalence $\Lambda(H)\to \Lambda(G)$. Thus, if
$\underline{X}$ is the orbifold represented by $G$, there is a well-defined orbifold $\Lambda(\underline{X})$
represented by $\Lambda(G)$. Following usage in algebraic geometry we call $\Lambda(\underline{X})$ the inertia
orbifold of $\underline{X}$.

For an orbifold groupoid of the form $K\rtimes M$ (see below 3.1) the space of loops is the ``Brylinski space''
\[S=\{(k,x)\,|\,x\in M, k\in K, k x=x\}\]
with its natural action by $K$ ($l\cdot(k,x)=(l k l^{-1}, l x)$). The orbit space $S/K$ was already considered in
\cite{K}.

Consider an arbitrary proper \'etale groupoid $G$ and a $G$-sheaf $B$. Then for the derived functor of
$\beta:\Lambda(G)\to G$, the discussion in 6.3 gives $R^{i} \beta_{*}(B)_{x}=H^{i}(G_{x}, B_{x})$ where
$G_{x}=\beta^{-1}(x)$ is viewed as a discrete set. Thus $R^{0}\beta_{*}=\beta_{*}$ and $R^{i}\beta_{*}=0$ for
$i>0$. So for the composite $q \beta:\Lambda G\to |G|$ we find
\[R^{i}(q\beta)_{*}(B)_{x}=R^{i}q_{*}(\beta_{*}(B))_{x}=H^{i}(G_{x}\rtimes G_{x},B_{x}).\]
Here $G_{x}$ acts on itself from the right by conjugation. Since the discrete groupoid $G_{x}\rtimes G_{x}$ is
equivalent to the sum of the centralizer subgroups $Z(g)$ indexed by conjugacy classes $(g)$ of elements $g\in
G_{x}$, we find
\[R^{i}(q\beta)_{*}(B)_{x}=\prod_{(g)}H^{i}(Z(g),B_{x}),\]
which fits into a Leray spectral sequence
\[H^{j}(|G|,R^{i}(q\beta)_{*}(B))\Rightarrow H^{i+j}(\Lambda(G),B).\]
In particular, if $B$ is the constant sheaf $\CC$ of complex numbers, one finds that
$R^{i}(q\beta)_{*}(\CC)_{x}=0$ for $i>0$ while $R^{0}(q\beta)_{*}(\CC)_{x}=\prod_{(g)}\CC=Class(G_{x},\CC)$, the
set of class functions on $\CC$, which is the same as the complex representation ring $R_{\CC}(G_{x})$. Thus
\[H^{i}(\Lambda(G),\CC)=H^{i}(|G|, \underline{R}_{\CC}),\]
where $\underline{R}_{\CC}$ is the representation ring sheaf with stalk $R_{\CC}(G_{x})$ at $x$. One might call
this the \textit{Bredon cohomology} of the orbifold represented by $G$. (Indeed, if $G$ represents a reduced
orbifold, then $G$ is Morita equivalent to a translation groupoid $K\rtimes M$ by Theorem 3.7, and it is known in
this case that $H^{i}(M/K,\underline{R}_{\CC})$ is isomorphic to the Bredon cohomology of $M$ with coefficients in
the representation ring system \cite{Ho}.)

If $|G|$ is compact, there is a Chern character isomorphism $(\nu=0,1)$
\[K^{\nu}(G)\otimes\CC\overset{\sim}{\to}\prod_{i}H^{2i+\nu}(|G|,\underline{R}_{\CC})\]
which is ``locally'' the one of \cite{Sl}; see \cite{AR}. This Chern character factors naturally as a composition
\[K^{\nu}(G)\otimes\CC\to H P^{\nu}(C^{\infty}_{c}(G))\to \prod_{i}H^{2i+\nu}(|G|,\underline{R}_{\CC}).\]
The first is the non-commutative Chern character into the periodic cyclic cohomology of the convolution algebra
$C^{\infty}_{c}(G)$ of $G$ (\cite{BC}). This is the algebra of compactly supported smooth functions
$\alpha:G_{1}\to \CC$, with product defined exactly as for the group ring of a finite group:
\[(\alpha\cdot\beta)(g)=\sum_{g=h k}\alpha(h)\beta(k)\]
where the sum is over all ways of writing an arrow $g$ in $G$ as a composite $g=h k$ in $G$. The second map is
essentially the isomorphism occurring in a somewhat different form in \cite{BN,Cr,CM3}.

\noindent
Ieke Moerdijk\\
Mathematical Institute,\\
University of Utrecht,\\
PO Box 80.010,\\
3508 TA \ Utrecht,\\
The Netherlands\\
Email: moerdijk@math.uu.nl

\end{document}